\documentclass{article}
\usepackage[utf8]{inputenc}
\usepackage{authblk}
\usepackage{setspace}
\usepackage[margin=1.25in]{geometry}

\usepackage{subcaption}
\usepackage{amsmath}
\usepackage{lineno}
\usepackage{hyperref}

\usepackage{graphicx}
\usepackage{wrapfig}
\usepackage{booktabs,tabularx}

\usepackage[style=nejm, 
citestyle=numeric-comp,
sorting=none,
backend=bibtex]{biblatex}
\title{
    On AI and Progress in Mathematics\\[0.2em]
    \large\textit{Dreaming a Mathematical Universal Machine}
}
\author[1*]{Vincent Runge}

\affil[1]{Université Paris-Saclay, CNRS, Univ Evry, Laboratoire de Mathématiques et Modélisation d'Evry, 91037, Evry-Courcouronnes, France.}
\affil[*]{vincent.runge@univ-evry.fr}

\date{}

\begin{document}

\maketitle

\begin{abstract}
In the autumn of 2026, mathematics is facing what may be the deepest crisis in its history. The entire ecosystem through which mathematical knowledge is produced, validated, communicated, and assessed is beginning to break down. Mathematics departments, those usually quiet and civilized places, are suddenly in turmoil -- the unlikely theaters of a new Battle of the Ancients and the Moderns. The source of the crisis is clear: AI systems have become better than even the best human mathematicians at proving and validating mathematical results. 

This document proposes a metaphor for mathematical activity: we once thought of ourselves as mountain climbers; soon, we may become miners digging in the mathematical mine, either collaborating or competing with mathematical AI machines. Publications would then be seen primarily as a means of helping humans understand mathematics, while the measurement of mathematical progress would be entrusted to a Mathematical Universal Machine. In nontechnical terms, we explore what this brave new world of mathematics might look like and how we might collectively shape its future.
\end{abstract}

\vspace{3cm}

\begin{flushright}
\begin{minipage}{0.62\textwidth}
\raggedleft
\emph{``La liberté est indivisible : on ne peut en retrancher une partie
sans la tuer tout entière.''}

\medskip
--- Mikhaïl Bakounine. 1895
\end{minipage}
\end{flushright}

\newpage 
\section{The Human Adventure}
\label{sec:adventure}

Freedom and creativity are among the values that mathematicians cherish most. Year after year, we patiently build inner worlds in which abstract concepts gradually become familiar and, strangely, almost concrete. \textbf{Step by step, we become explorers of our own small planets}, learning every path, every landmark, and every hidden corner. We take pride in becoming a reference, even a modest one, on a particular subject. Within this small territory, we gain the capacity to judge new ideas.
\begin{wrapfigure}{r}{0.30\textwidth}
    \centering
\includegraphics[width=0.30\textwidth]{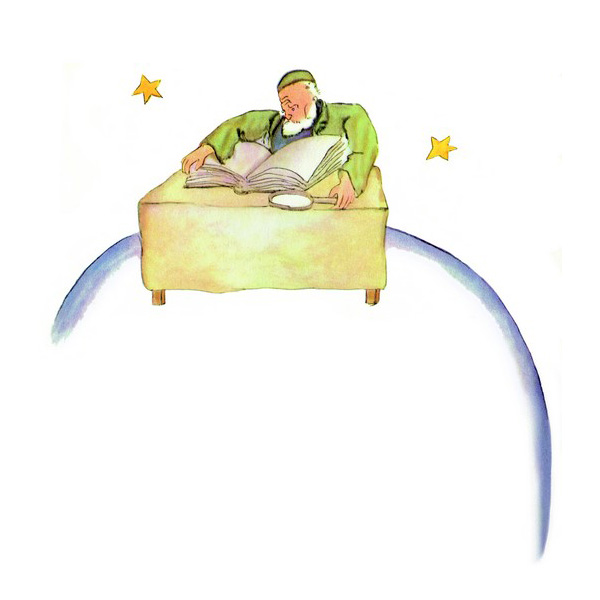}
    \caption*{The Geographer, from Antoine de Saint-Exupéry’s {\it Le Petit Prince} (1943).}
    \label{fig:geographer}
\end{wrapfigure}
We share our enthusiasm for what is genuinely new and patiently explain what is not. We encourage new generations of students to risk an adventure in the beautiful landscapes of mathematics. Yet with AI, we risk becoming like the Geographer in Saint-Exupéry's famous book, who never leaves his desk and cannot even say whether his own planet has oceans or mountains: he waits for explorers to tell him. Tomorrow, the explorers may be machines, restricting our freedom and drying up our creativity.

We also love sharing this knowledge with our colleagues, and much of what we share consists of surprisingly intimate thoughts: intuitions that are difficult to formalize, doubts, half-formed analogies, and memories of failures. Some of these colleagues become close friends. We drink a lot of coffee\footnote{We know the famous quote, ``A mathematician is a device for turning coffee into theorems,'' attributed to Alfréd Rényi but often mistakenly ascribed to Paul Erdős.} and stand in front of blackboards for hours. We draw pictures and invent metaphors, trying to make one another see what we have begun to see. Sometimes we meet the giants: mathematicians whose names we have known for years, from papers that once seemed impossibly difficult. We thank them; we admire them. And then, sometimes, one of them reads something we have written, and for a brief moment we see a small light of recognition in their eyes. \textbf{We inherit stories as well as theorems.} Late in the evening, at a conference dinner, an older mathematician tells us an anecdote about Erdős, Grothendieck, or Arnold. Such stories are not in the books, yet they connect us to the deep past of mathematics. They remind us that mathematics is an immense construction whose builders, even after they are gone, somehow remain present among us.

Every mathematician is, in this sense, the curator of an old and peculiar museum to which they hold a key. When they grant us the privilege, they take us on a private tour. They show us not only what is displayed in the galleries but also why it matters, how it came to be there, which doors lead nowhere, and which forgotten passage might still be worth opening. \textbf{Books, papers, and Wikipedia contain mathematics, but they contain it cold.} The living experience of mathematics is something else: it lies in learning from others and sharing what we have learned, in the transmission of intuition, taste, stories, enthusiasm, and ways of seeing. All of this, and much more, is the human adventure of mathematics. With a Machine writing our theorems, would we still care, connect, and engage in this adventure?

\newpage 

\begin{wrapfigure}{l}{0.38\textwidth}
    \centering \includegraphics[width=0.37\textwidth]{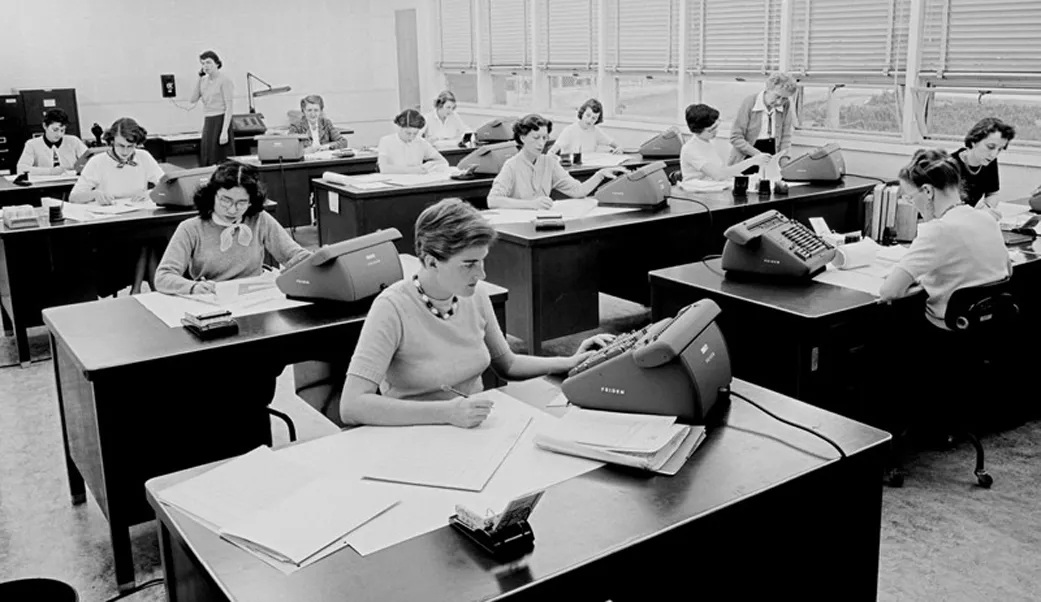}
    \caption*{A talented team of women known as ``computers'' at the Jet Propulsion Laboratory.}
    \label{fig:JPL}
\end{wrapfigure} 
Mathematical activity is very hard. Generation after generation, we have patiently built the coral reef of mathematics: concept upon concept, definition upon definition, proof upon proof. We fill the gaps between fields, build bridges between distant ideas, and slowly transform isolated results into coherent bodies of knowledge. AI has changed this game. \textbf{We were trained to prove, and suddenly proof itself risks becoming obsolete as a distinctively human skill.} It is as if we had spent years learning to multiply large numbers by hand, only for someone to place a calculator on the table. This analogy is less metaphorical than it first appears. For centuries, computation was itself a human activity. During the first half of the twentieth century, women were employed as ``human computers,’’ performing calculations by hand until electronic computers gradually took over this work. Some of these women went on to become among the earliest computer programmers. Something similar may now be happening to mathematicians in our proof-oriented activity: what machines first did to computation, they may now be beginning to do to proof. We are moving from an era dominated by \emph{Human Understanding} toward one in which we must contend with \emph{Machine Understanding}. New proofs, connections, and interpolations between mathematical fields can now be produced at a rate that may soon exceed anything the mathematical community has experienced. The scarcity on which much of our mathematical culture was built is beginning to disappear. \textbf{When proofs become abundant, where is the new scarcity?}\\
 
\begin{wrapfigure}{r}{0.35\textwidth}
    \centering
\includegraphics[width=0.17\textwidth]{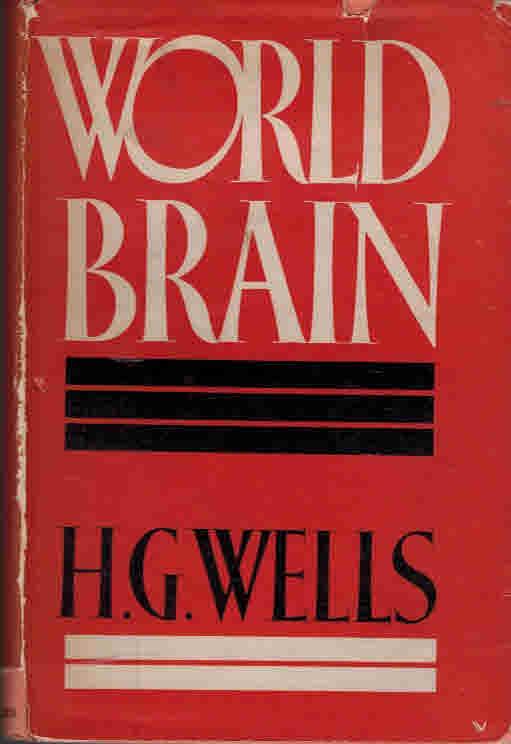}
\includegraphics[width=0.17\textwidth]{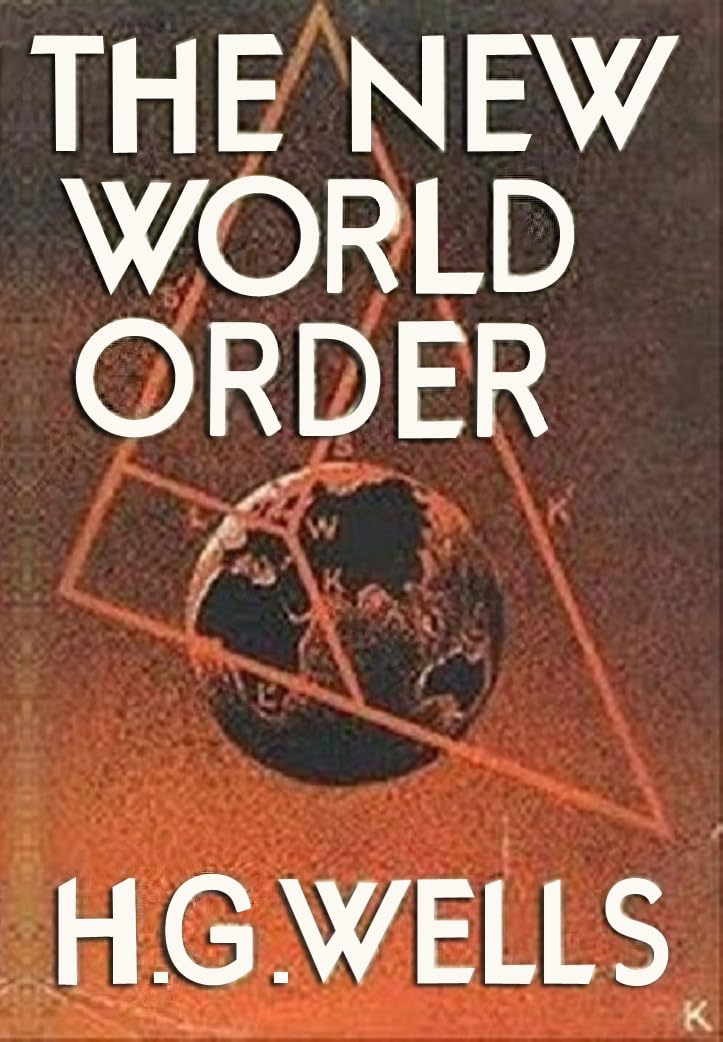}
    \caption*{H.G.~Wells, \emph{World Brain} (1938) and \emph{The New World Order} (1940).}
        \label{fig:wells}
\end{wrapfigure}

None of this introduction will be news to mathematicians.  It simply spells out why the present transformation feels so threatening. \textbf{We fear breaking with a centuries-old tradition of intellectual life} and entering an obscure and uncertain world in which human mathematical understanding is gradually left behind. We want to preserve what we value most: intellectual independence, freedom, universality, openness, and the pleasure of understanding mathematics for ourselves.

But there is another danger. Mathematics rests on unusually strong norms of trust, integrity, and intellectual honesty. \textbf{AI risks releasing one of the worst genies from the bottle: suspicion.} A divide between the AI-free mathematician and the AI-assisted mathematician could fracture the community. We therefore need an environment in which human mathematical activity can be both preserved and transformed. The old system of papers, peer review, journals, citations, and prestige was not perfect; faced with an abundance of machine-generated mathematics, it may become increasingly inadequate. \textbf{Our proposal is to partially decouple the advancement of human understanding from the measurement
of mathematical progress}: to protect the former while developing new tools to measure and organize the latter.\\

This requires a new institution for this new world order: a Mathematical Universal Machine (\textsc{MUM}), a shared AI infrastructure available to the entire mathematical community. It would extend H.G.~Wells's vision of a \emph{World Brain}~\cite{wells1938} (see Appendix~\ref{app:shared} for an AI-generated historical overview!). \textbf{This Machine would also provide quantitative measures of mathematical progress} that complement the traditional peer-review system. Designing the Machine and measuring its behavior is, in itself, an enormous mathematical task. We present below what such a system might look like.

\section{The Mathematical Universal Machine}
\label{sec:machine-update}

\subsection{The mine and the miner: maps, tools, brain and budget}

In this new world of AI -ssisted mathematics, we propose the metaphor of the {\bf mathematical mine}. The mine is the space of mathematical truth. Mathematical results are diamonds hidden within it; proofs are tunnels leading to them. The mathematical AI system is the {\bf miner}. Its purpose is to navigate this space by making its own decisions with some degree of independence in order to complete the task. The presentation is restricted to very generic terms. An AI-assisted search led us to much more precise ideas that anyone can easily reproduce. 

We consider a miner composed of {\bf three distinct components}:
\begin{enumerate}
\item The miner holds {\bf maps} leading to known diamonds. They correspond to research papers, databases of proven results, definitions, concepts, graphs of structured dependencies between theorems, etc. He can use them at any moment to improve his search capacity.
\item The miner can also use external {\bf tools}, such as Lean for formal verification, Python and computer algebra for experimentation. These tools make his work in the mine easier, helping him explore, verify, test, and orient himself.
\item The miner has a {\bf brain}: the core of the AI system, a large language model with billions of learned parameters. It enables him to work autonomously, searching for proofs, testing conjectures, using external tools, and deciding what to explore next with its own exploration strategy. It embodies his internalized mathematical intuition, representations, heuristics, and ability to decide where and how to dig. Importantly, the miner also has a memory and remembers the actions he has already taken, allowing him to learn from failed attempts and avoid repeatedly exploring the same paths.
\end{enumerate}

External to the miner's structure, the Machine has finite patience, represented by a finite computational {\bf budget}. His energy is limited by the human-designed system, whose computational resources must ultimately be shared among all mathematicians. The limitation may also be deliberate: we may choose to prevent excessively long explorations that could carry the miner too far into unknown territory. We may also choose to build a minimal-energy system, whose ratio of mathematical capacity to energy consumption is the new objective of mathematical research.

With all these components fixed, the AI system can perform some {\bf work} assigned to it. At the same time, letting such a system freely explore mathematics leads to what we call a {\bf convexification of the mine}: the Machine progressively fills in the mathematical territory accessible from what is already known. New results are discovered, new theorems are proved, new definitions emerge, simpler proofs are found, and connections between previously separated regions are uncovered. All of this new mathematics is made possible thanks to the Machine.

\textbf{We hypothesize, however, that human intelligence and imagination will not be entirely subsumed by the Machine.} Some mathematical work may still lie outside this convexification: new ways of presenting known mathematics, unexpected definitions or viewpoints, or genuinely original research that the Machine would not have discovered on its own. Once integrated into the Machine, such contributions enlarge its accessible territory and initiate a new process of convexification. This perspective suggests two distinct and complementary modes of mathematical research:

\begin{enumerate}
\item \textbf{The mathematician as a supervisor.}
The mathematician explores mathematics through the AI system, supervising the miner in the mine. The user commands the AI brain and constrains his work to his own expertise as a mathematician, telling it what to explore and proposing search directions and ideas. While the Machine seeks new results with human help and feedback, it constantly learns how to better use its capacity among tools, maps, memory, computations, etc. In practice, this may develop into an art of prompting the AI system, involving complex sequences of instructions, feedback, redirections, and interventions as the exploration unfolds.

\item \textbf{The mathematician as a miner.}
The mathematician works outside the Machine, exploring the mine himself and seeking mathematics beyond its current reach. The resulting papers provide new maps which, through \textbf{a Machine Update}, may be internalized into the miner's brain and change how it subsequently explores mathematics. In this mode, the mathematician deliberately limits his reliance on the Machine: following it too closely risks remaining within territory the Machine has already explored rather than venturing beyond its current convexification.

\end{enumerate}

These two modes represent complementary approaches to mathematical research: either supervising and guiding the explorer of the mine or remaining a miner oneself and pursuing mathematics through a more traditional form of exploration. To support both forms of research, we need the introduction of new metrics for evaluating mathematical production and quantifying the contributions of mathematicians in each role.

\newpage 

\subsection{The Machine horizon}

\begin{wrapfigure}{r}{0.45\textwidth}
    \centering
\includegraphics[width=0.45\textwidth]{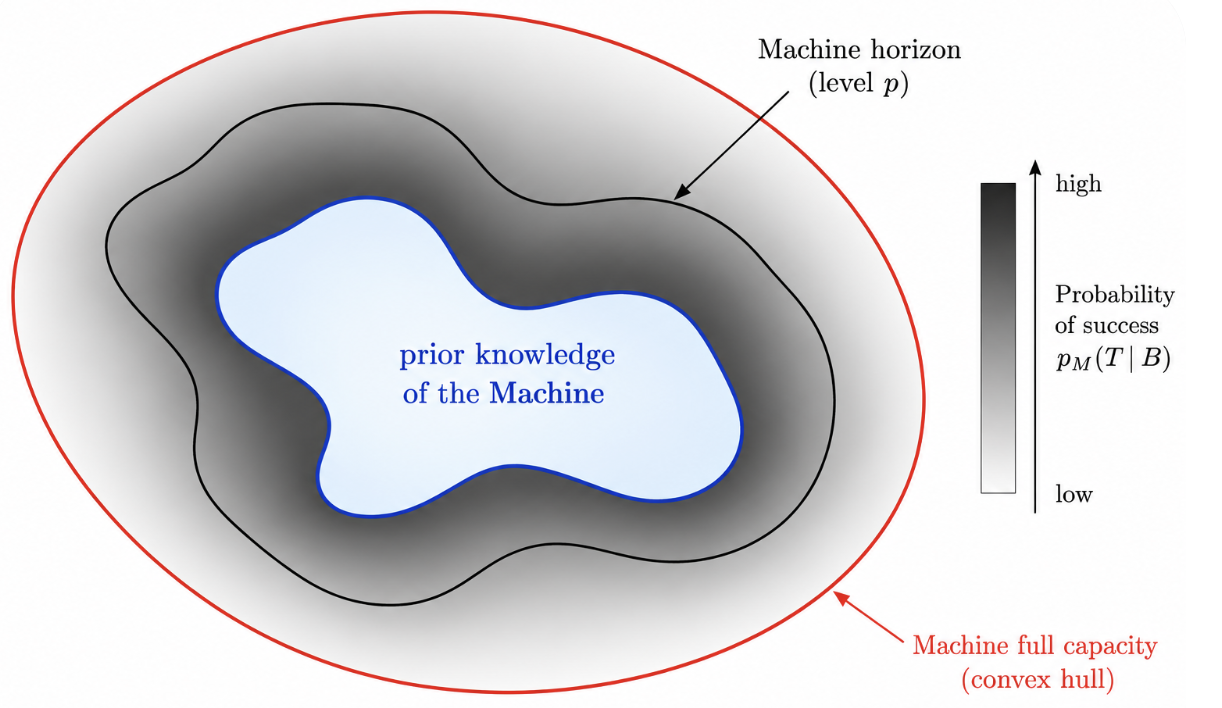}
    \caption*{}
    \label{fig:MUM}
\end{wrapfigure}

We formally represent the Mathematical Universal Machine, or {\bf miner}, by
$M=(\theta,\mathcal K,\mathcal T)$, where $\theta$ denotes its {\bf brain}, $\mathcal K$ its mathematical {\bf maps}, and $\mathcal T$ its external {\bf tools}. The miner operates under a finite computational budget $B$, measured in Mathematical Research Compute Units (MRCU): a quantity that requires a precise definition. For a target statement $T$, we define the probability of success as
$$
p_M(T\mid B)
= \Pr(M \text{ proves } T \text{ within budget } B)\,,
$$
where $T$ is known to humans but has been excluded from the Machine’s prior knowledge. For a given probability $p\in(0,1)$, the {\bf Machine horizon} is the set of mathematical results that the miner can reach within its computational budget:
$$
\mathcal H_M(B;p) =
\left\{T:p_M(T\mid B)\geq p\right\}.
$$
This horizon represents the mathematical territory accessible to the Machine: beyond it lie truths that the miner cannot reliably discover with its current brain, maps, tools, and computational resources. (On the history of this idea, see Appendix~\ref{app:reasoning}.)

\begin{wrapfigure}{r}{0.20\textwidth}
    \centering
\includegraphics[width=0.20\textwidth]{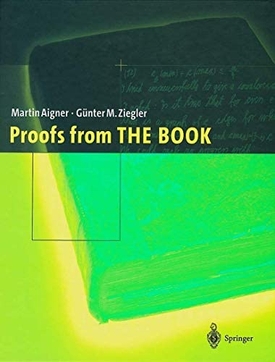}
    \caption*{Proofs from THE BOOK, first edition (1998)}
    \label{fig:book}
\end{wrapfigure}
Extending the Machine horizon $\mathcal H_M(B;p)$ echoes Paul Erdős’s vision of {\it The Book}, an imaginary book in which God keeps the most beautiful proof of every mathematical theorem\footnote{As Paul Erdős famously put it, ``You don’t have to believe in God, but you should believe in The Book.’’}.  This vision inspired Martin Aigner and Günter M. Ziegler’s {\it Proofs from THE BOOK}, a collection celebrating particularly elegant mathematical proofs. Our proposal pursues a related goal: rather than seeking the most beautiful proof of each theorem, we seek to make mathematical truths increasingly accessible by improving the miner’s ability to discover them. Just as a proof from {\it The Book} reveals an unexpectedly elegant path to a mathematical truth, a more capable Machine may uncover shorter paths to entire regions of the mathematical mine.\\

Improving the Machine requires new metrics to assess its progress. Today, the work of mathematicians is evaluated through a long human-based peer-review process. As said earlier, publication could serve the purpose of advancing human understanding rather than measuring mathematical progress. The Mathematical Universal Machine would provide quantitative measures of mathematical progress and help attribute individual and collective contributions to mathematical discoveries.

The Machine would be updated at fixed intervals by incorporating new publications, leading not only to new mathematical maps but also to an update of the miner’s brain parameters, $\theta_t \to \theta_{t+1}$. After each {\bf Machine Update}, its mathematical capabilities would be reassessed, notably by measuring changes in the probability of proving a theorem, $p_{M_t}(T\mid B)$, and the expansion of its horizon, $\mathcal H_{M_t}(B;p)$. We could also measure reductions in the computational effort required to discover mathematical results, thereby estimating how much a contribution improves the Machine’s ability to explore the mine. \textbf{Such measurements could provide a new basis for recognizing the contributions of mathematicians}, whether they discover new mathematics themselves or improve the Machine’s exploration through human supervision.\\
 
This vision of a new organization of mathematics, in which industrial and handcrafted mathematical work coexist, remains speculative, as does the possibility of developing a system that reconciles the two. Its realization will depend on our ability to develop mathematical AI systems whose capabilities and limitations are better understood.

\section{The Political Project}
\label{sec:political}

Mathematicians often share progressive views on society, science, and technology. As scientists and citizens, many of us also care deeply about ecology and climate change. The current trajectory of AI is not what we dreamed of. AI consumes considerable amounts of energy, water, land, and physical resources. Its use also raises serious ethical concerns when associated with surveillance, social control, or warfare. AI for mathematics is currently concentrated in the hands of a few private companies, mostly based in the United States, whose most powerful models are accessible only to selected collaborators. Meanwhile, the negative physical impacts of using these systems for mathematics remain poorly measured and documented.  Most mathematicians find the present state of affairs difficult to accept.\\

We must not resign ourselves to the status quo. We can organize, act collectively, and propose a better model. {\bf We can build a universal AI infrastructure for mathematics, developed in partnership with an international institution such as UNESCO} and governed by mathematicians through a dedicated international organization (for instance, the International Mathematical Union). Supported by states, universities, and philanthropic institutions, it would belong to the global scientific commons rather than remain a proprietary asset controlled by a handful of private actors. This Mathematical Universal Machine (\textsc{MUM}) would be built around the values of the mathematical community. As a starting point, it could be governed by the following principles:

\begin{enumerate}
\item \textbf{Ecological accountability.} Its energy consumption, environmental impact, and material footprint would be independently assessed and published every year.
\item \textbf{Open architecture.} The architecture of the model, its weights, and the essential components of the system would be openly accessible.
\item \textbf{Universal and equal access.} Access would be open to the mathematical community -- potentially including students -- with an equal computational or token budget for every user.
\item \textbf{Privacy.} Users would retain full control over their data. Private research, unpublished results, and preprints submitted to the system would never be incorporated into the shared system without explicit consent.
\item \textbf{Transparent Machine Updates.} Once a year, the mathematical database and the miner's brain would be updated using newly published mathematics and eligible preprints. The computational importance of each new and older contribution would then be estimated.
\item \textbf{International scientific governance.} A joint mathematics--computer science laboratory, established within or in partnership with UNESCO, would be responsible for the Machine Update, maintenance, evaluation, and communication.

\end{enumerate}

Mathematicians often describe their research through the metaphor of climbing a mountain. Through an immense solitary effort, we reach a summit from which we contemplate a vast landscape of unexplored mathematical truths. As Terence Tao has suggested, we may now be witnessing the end of the romantic era of mathematics, the age of solitary heroes. \textbf{From the mountain, we descend into the mine:} an industrial landscape filled with machines, unexpected discoveries, and disruptive innovations.
According to a famous anecdote, when John von Neumann was asked what percentage of mathematics a person could aspire to understand, he replied: ``About 28 percent.'' With the Machine, we may recover something of the nineteenth-century ideal, when leading mathematicians could still aspire to a comprehensive understanding of mathematics as a whole.

In the mud of the mine, we have the power to collectively organize the search, preserve our ideals, and perhaps build the era of industrial mathematics ourselves. Amid this transition, one truth remains: we have been sent from the mountain into the mine. For progress, there is no cure\footnote{John von Neumann. ``Can We Survive Technology?'', Fortune, June 1955.}. Adaptation is not optional; preserving our values and some of our working habits is still possible.

\newpage 

\appendix

The AI suggested including these appendices in the document and generated them, with only minor manual corrections.

\section{Shared knowledge}
\label{app:shared}

Each generation of mathematicians has reinvented the way it shares what it knows, and with it the way credit is given. This history reads as the slow unbundling of the scientific journal (Table~\ref{tab:shared}). The \textsc{MUM} would be its next step, and the first in which the Machine itself is a contributor.

Before journals, letters carried mathematics: Marin Mersenne's correspondence earned him the name of ``post-box of Europe''. In 1665, the \emph{Journal des s\c{c}avans} (Paris, 5 January) and the \emph{Philosophical Transactions} (London, 6 March) bundled four functions we still rely on: registering who found what first, certifying results, making them known, and archiving them~\cite{roosendaal1997}. The \emph{Encyclop\'edie} of Diderot and d'Alembert (1751--1772)~\cite{encyclopedie} added a map, d'Alembert's tree of knowledge, and linked its articles by cross-references; when its privilege was revoked in 1759, the last volumes were completed clandestinely. From 1868, the \emph{Jahrbuch über die Fortschritte der Mathematik}, then Neugebauer's \emph{Zentralblatt} (1931) and \emph{Mathematical Reviews} (1940), mapped the whole literature, now freely accessible as zbMATH Open; databases of objects such as the OEIS (1964)~\cite{oeis} followed. They are the ancestors of the maps $\mathcal K$ of Section~\ref{sec:machine-update}.

In the 1930s, Bourbaki rewrote mathematics on a common architecture of structures under a collective pseudonym~\cite{bourbaki}, and became a monoculture: whoever designs a shared architecture shapes what counts as mathematics. H.G.~Wells's \emph{World Brain}~\cite{wells1938} called for continual revision and for copies spread around the world, two ideas that return as the annual Machine Updates and the open architecture of Section~\ref{sec:political}, without the centralism later readers criticized~\cite{rayward1999}.

The arXiv (1991)~\cite{ginsparg2011} separated dissemination from certification: preprints circulate at once and timestamps settle priority, leaving journals with certification alone. The decoupling proposed in Section~\ref{sec:adventure} splits certification in turn: proof assistants check correctness, the \textsc{MUM} would measure progress, and publications would serve human understanding. Polymath (2009)~\cite{gowers2009} made public the intimate side of research described in Section~\ref{sec:adventure} and found a new proof of the density Hales--Jewett theorem, signed D.H.J.~Polymath~\cite{polymath2012}; but no hiring committee can evaluate a share of a pseudonym. The Equational Theories Project (2024) settled all $22{,}028{,}942$ implications between $4{,}694$ laws of magmas by combining human proofs, automated provers and Lean, and documented who contributed what~\cite{etp2025}: a small-scale prototype of the \textsc{MUM}.

Credit, meanwhile, came to be counted, with citation indexes~\cite{garfield1955} and the h-index~\cite{hirsch2005}; in its report \emph{Citation Statistics}, the International Mathematical Union warned against ranking with such numbers~\cite{adler2009}. The \textsc{MUM} would not count citations but estimate the difference a work makes to what the Machine can reach (Section~\ref{sec:machine-update}), as a complement to human judgment. Mathematicians have also acted collectively when access to knowledge was threatened: the editors of \emph{Topology} resigned in 2006~\cite{jackson2007}, and thousands joined ``The Cost of Knowledge'' boycott in 2012~\cite{costofknowledge}. Building a public Machine (Section~\ref{sec:political}) is the same fight with a new tool.

\begin{table}[ht!]
\centering
\small
\begin{tabularx}{\textwidth}{@{}l l >{\raggedright\arraybackslash}X >{\raggedright\arraybackslash}X@{}}
\toprule
\textbf{Step} & \textbf{Date} & \textbf{Who contributes, and how} & \textbf{Lesson for the \textsc{MUM}}\\
\midrule
Journals & 1665-- & Editors and correspondents; printed issues & Registration, certification, dissemination and archiving, bundled together\\
\emph{Encyclop\'edie} & 1751--1772 & A society of men of letters; volumes linked by cross-references & Knowledge needs a map; sharing it is political\\
Reviews, zbMATH & 1868-- & Reviewers covering the whole literature & The first maps of the whole discipline\\
Bourbaki & 1935-- & A secret collective; one unified treatise & A common architecture, and the risk of monoculture\\
\emph{World Brain} & 1936--1938 & A world editorial staff; microfilm copies & Continual revision and replication, without centralism\\
Citation indices & 1964-- & Every citing author, counted & Credit by counting, and its misuse\\
arXiv & 1991-- & Anyone with a preprint; a public server & Dissemination separated from certification\\
Polymath & 2009-- & An open crowd; blog and wiki & Collective work succeeds, but credit dissolves\\
Equational Theories & 2024--2025 & Humans, automated provers and Lean & Human and machine work combined, with documented contributions\\
\textsc{MUM} & --- & Humans and the Machine; weights, tools and maps & Credit can no longer be read off an author list: it must be measured\\
\bottomrule
\end{tabularx}
\caption{Sharing mathematical knowledge, from the first journals to the \textsc{MUM}.}
\label{tab:shared}
\end{table}

\newpage

\section{Mechanizing mathematical reasoning}
\label{app:reasoning}

The second tradition behind the \textsc{MUM} asks whether reasoning, like computation before it (Section~\ref{sec:adventure}), can be mechanized. Its answer is the model of Section~\ref{sec:machine-update}: not a procedure deciding every statement, but a search with a budget and a probability of success. Read in order, it supplies the components of the miner $M=(\theta,\mathcal K,\mathcal T)$ and his budget $B$ (Table~\ref{tab:reasoning}).

\begin{table}[ht!]
\centering
\small
\begin{tabularx}{\textwidth}{@{}l l >{\raggedright\arraybackslash}X >{\raggedright\arraybackslash}X@{}}
\toprule
\textbf{Step} & \textbf{Date} & \textbf{What it established} & \textbf{Component of the \textsc{MUM}}\\
\midrule
Leibniz, \emph{calculemus} & c.\ 1685 & Reasoning as calculation; disputes settled by checking & Formal statements; an answer to suspicion\\
Frege; \emph{Principia} & 1879; 1910--1913 & Mathematics in a fully formal language & The language of the maps, and its cost\\
Gödel; Turing & 1931; 1936 & Incompleteness; undecidability; the universal machine & Truths out of reach; the budget $B$\\
Gödel's letter; Cook & 1956; 1971 & The cost of finding a proof (P versus NP) & The horizon $\mathcal H_M(B;p)$\\
Logic Theorist; EQP & 1956; 1996 & Machines find proofs, then settle an open problem & Search, and a first question of machine credit\\
Automath, Mizar, Lean & 1967--2015 & Proof checkers & The tools $\mathcal T$\\
Four Color; Kepler & 1976--2014 & Computer proofs, then formal certification & Suspicion answered by verification\\
QED, mathlib & 1994; 2017-- & A shared, verified library & The maps $\mathcal K$ and their governance\\
Learned provers & 2020--2025 & Search guided by learned intuition & The brain $\theta$, learned from the maps\\
\bottomrule
\end{tabularx}
\caption{Mechanizing mathematical reasoning, from Leibniz to the \textsc{MUM}.}
\label{tab:reasoning}
\end{table}

Leibniz hoped that disputes would be settled by calculation, the parties saying to each other \emph{calculemus}, ``let us calculate''~\cite[200]{leibniz}: an early answer to the suspicion of Section~\ref{sec:adventure}, since a checked proof requires no trust. Frege, then Whitehead and Russell, built the formal language~\cite{frege1879,principia} and revealed its cost: $1+1=2$ follows only from proposition $\ast$54.43 of \emph{Principia Mathematica}. Hilbert asked for a procedure deciding every statement~\cite{hilbert1928}. Gödel showed that any consistent, effectively axiomatized theory containing arithmetic leaves some truths unprovable~\cite{godel1931}: from given maps, some diamonds stay out of reach at any budget. Turing showed that no algorithm decides provability~\cite{turing1936}. Proofs can still be listed one by one, so an unlimited search eventually finds every provable statement, but nobody can tell in advance when, or whether, a given statement will appear: this is the limit of the convexification of Section~\ref{sec:machine-update}, and the reason why a search must be cut off. In the same paper Turing described the universal machine; the \textsc{MUM} is universal in Turing's sense and in that of Article~27 of the Universal Declaration of Human Rights, the right of everyone ``to share in scientific advancement and its benefits''~\cite{udhr1948}, which grounds the equal access of Section~\ref{sec:political}.

Undecidability left open the question that matters in practice: the cost of finding a proof. In 1956, Gödel wrote to von Neumann, then gravely ill: a machine can easily decide whether a statement has a proof of length $n$, but how fast must its running time grow with $n$? If only linearly or quadratically, ``the mental work of a mathematician concerning Yes-or-No questions could be completely replaced by a machine''~\cite{sipser1992}. The question resurfaced as the P versus NP problem, first stated for theorem-proving procedures~\cite{cook1971}, and is still open: checking a proof is easy, finding one seems hard. This asymmetry organizes the \textsc{MUM}: verification is cheap and certain, the work of the tools $\mathcal T$; search is costly and uncertain, the work of the brain $\theta$ under the budget $B$. The horizon $\mathcal H_M(B;p)$ asks Gödel's question empirically: not how the cost grows for the worst statement, but how far a given miner reaches with his budget.

Machines then began to search and to check. The Logic Theorist of Newell, Shaw and Simon (1956)~\cite{newell1956} proved 38 of the first 52 theorems of \emph{Principia}, one of them more elegantly than Whitehead and Russell, an early echo of \emph{The Book} (Section~\ref{sec:machine-update}); a paper listing the program as a co-author was reportedly refused. In 1996, the prover EQP settled the Robbins problem, open for sixty years~\cite{mccune1997}: a machine had reached a diamond that humans had not. Meanwhile, computation entered proofs, and with it suspicion: no referee could redo by hand the computer check behind the Four Color Theorem (1976), and the referees of the Kepler conjecture were only ``99\% certain''. Proof checkers, from Automath and Mizar to Coq and Lean~\cite{debruijn1970,demoura2015}, removed the doubt with formal proofs of both theorems (2005 and 2014)~\cite{gonthier2008,hales2017}, moving trust from referees to verification. The QED Manifesto (1994) called for a single verified library of mathematics~\cite{qed1994} and stalled~\cite{wiedijk2007}; mathlib has since built one in the open~\cite{mathlib2020}, and the Liquid Tensor Experiment (2022) took it to current research~\cite{scholze2022}.

The brain came last. The heuristics of the Logic Theorist were written by hand; those of today's provers are learned from the maps. In 2020, a language model trained on formal mathematics guided proof search and contributed proofs to the Metamath library~\cite{polu2020}; in 2025, AI systems reached gold-medal level at the International Mathematical Olympiad~\cite{deepmind2025}. Because the brain is learned from the maps, new maps change the brain: this is the Machine Update $\theta_t\to\theta_{t+1}$ of Section~\ref{sec:machine-update}. Such provers are already scored by the share of benchmark statements they prove within a fixed number of attempts~\cite{zheng2022}: the probability $p_M(T\mid B)$ in all but name. The \textsc{MUM} would add scale and openness: the whole of mathematics as benchmark, and an institution that keeps maps, tools, brain and budget open and shared (Section~\ref{sec:political}).

Turing had named the open question. In 1939 he described mathematical reasoning as ``the exercise of a combination of two faculties, which we may call intuition and ingenuity'', recalling that ``in pre-Gödel times'' some had hoped to replace all intuitive judgments by formal rules~\cite{turing1939}. Gödel's letter asked whether ingenuity could be mechanized efficiently; learned provers suggest that part of intuition can be learned. Whether all of it can is the hypothesis of Section~\ref{sec:machine-update}, which the \textsc{MUM} is designed to test.

\printbibliography

\end{document}